\documentclass[11pt,english]{article}
\usepackage{amsfonts,amsmath,amsthm,amscd,amssymb,latexsym}

    \usepackage[T2A]{fontenc}
    \usepackage[cp1251]{inputenc}
    \usepackage[ukrainian]{babel}
    \usepackage{graphicx}

\begin{document}

\title{ПРО ОБЕРНЕНУ НЕРІВНІСТЬ ПОЛЕЦЬКОГО З КОТАНГЕНСАЛЬНОЮ ДИЛАТАЦІЄЮ}

\author{Євген О.\,Севостьянов, Валерій А.\,Таргонський}
%\shorttitle{Short paper title for the headers}
%\shortauthor{F. Author, S. Author}

%Change theorem environments according to your needs...
\theoremstyle{plain}
\newtheorem{theorem}{Теорема}[section]
\newtheorem{lemma}{Лема}[section]
\newtheorem{proposition}{Твердження}[section]
\newtheorem{corollary}{Наслідок}[section]
\theoremstyle{definition}

\newtheorem{example}{Приклад}[section]
\newtheorem{remark}{Зауваження}[section]
\newcommand{\keywords}{\textbf{Key words.  }\medskip}
\newcommand{\subjclass}{\textbf{MSC 2000. }\medskip}
\renewcommand{\abstract}{\textbf{Анотація.  }\medskip}
\numberwithin{equation}{section}

\setcounter{section}{0}
\renewcommand{\thesection}{\arabic{section}}
\newcounter{unDef}[section]
\def\theunDef{\thesection.\arabic{unDef}}
\newenvironment{definition}{\refstepcounter{unDef}\trivlist
\item[\hskip \labelsep{\bf Определение \theunDef.}]}%
{\endtrivlist}

\renewcommand{\figurename}{Мал.}

\maketitle

\begin{abstract}
Статтю присвячено встановленню спотворення модуля сімей кривих в
широких класах відображень, які допускають наявність точок
розгалуження. Зокрема, для відображень, які є диференційовними майже
скрізь, мають $N$- та $N^{\,-1}$-властивості Лузіна та є абсолютно
неперервними на майже всіх кривих, отримана обернена нерівність
Полецького з так званою котангенсальною дилатацією. Показано, що для
відображень, що мають обернені, ця дилатація співпадає з так званою
тангенсальною дилатацією оберненого відображення. Крім того,
обґрунтовано, що котангенсальная дилатація є меншою за максимальну,
зокрема, може бути меншою одиниці на множині додатної міри.
\end{abstract}

\medskip
{\bf Ключові слова:} модулі сімей кривих і поверхонь, обернена
нерівність Полецького, відображення зі скінченним та обмеженим
спотворенням

\medskip
\centerline{\bf English title: On the inverse Poletsky inequality
with cotangent dilatation}

\medskip
{{\bf English abstract.} {\bf On the inverse Poletsky inequality
with cotangent dilatation.} The article is devoted to establishing
the distortion of the modulus of families of paths in wide classes
of mappings that admit branch points. In particular, for mappings
that are differentiable almost everywhere and have $N$- and $N^{\,-
1}$-Luzin properties and are absolutely continuous on almost all
paths, we obtained the inverse Poletsky inequality with the
so-called cotangent dilatation. We have proved that, for inverse
mappings, this dilatation coincides with the so-called tangential
dilatation of the corresponding inverse mapping. In addition, we
have proved that cotangent dilatation is less than the outher or
inner dilatation, in particular, may be less than one on the set of
positive Lebesgue measure. }



\section{Вступ} Як відомо, оцінки модуля сімей кривих при
відображеннях відіграють ключову роль при їх вивченні, див., напр.,
\cite{Cr}, \cite{MRV$_1$}-\cite{MRV$_3$},\cite{MRSY}, \cite{Re} і
\cite{Ri}. Зокрема, вкажемо на публікації, в яких такі оцінки
використовують так звану тангенсальну (дотичну) дилатацію
відображення, див., напр., \cite{GG}, \cite{RSY}, \cite{RSSY} і
\cite{SalSev}. Тангенсальна дилатація була впроваджена з метою
узагальнення добре відомих внутрішньої і зовнішньої дилатації, а її
використання є важливим з точки зору вивчення їх локальної поведінки
і застосувань до рівняння Бельтрамі, див., напр., \cite{RSY}. Слід
зауважити, що на даний час відомі лише верхні оцінки спотворення
модуля сімей кривих з використанням тангенсальних дилатацій. Мета
даного рукопису -- довести аналогічні результати для сімей кривих в
прообразі при відображенні. Зокрема, за допомогою них можна
встановити порядок спотворення відстані при ньому, у тому числі,
логарифмічну неперервність за Гельдером, див. \cite{SSD} і
\cite{SevSkv}. Окремо буде розглянуто питання щодо отримання як
прямих, так і обернених оцінок спотворення модуля сімей кривих
плоских відображень, які використовують тангенсальну дилатацію і її
<<обернений>> аналог.

\medskip
Нагадаємо деякі означення, необхідні нам для формулювання одного з
основних результатів. Нехай $I=[a,b].$ Для спрямлюваної кривої
$\gamma:I\rightarrow {\Bbb R}^n$ визначимо функцію довжини
$s_{\gamma}(t)$ за наступним правилом: $s_{\gamma}(t)=S\left(\gamma,
[a,t]\right),$ де $S(\gamma, [a,t])$ позначає довжину кривої
$\gamma|_{[a, t]}.$ Для множини $B\subset I$ символом $S(\gamma, B)$
позначимо міру множини значень, що приймає функція $s_{\gamma}(t)$
на множині $B.$

\medskip
\begin{remark}\label{rem15}
Нехай $\gamma:[a, b]\rightarrow {\Bbb R}^n$ спрямлювана крива,
$t_0\in (a, b),$  $l_{\gamma}(t)$ позначає довжину підкривої $\gamma
|_{[t_0,t]}$ кривої $\gamma$ при $t>t_0,$ $t\in (a, b),$ і
$l_{\gamma}(t)$ дорівнює $-S(\gamma, [t,t_0])$ при $t<t_0,$ $t\in
(a, b).$ Зауважимо, що властивості функції $L_{\gamma, f}$ між
натуральними параметрами $l_{\gamma}(t)$ і
$l_{\widetilde{\gamma}}(t)$ кривих $\gamma$ і $\widetilde{\gamma}$
таких, що $\widetilde{\gamma}=f\circ\gamma,$ істотно не залежать від
обирання $t_0\in (a, b).$ У випадку замкненої кривої $\gamma$ ми
вважатимемо, що $t_0=a,$ оскільки при заданному $t_0\in (a, b)$
виконано рівність $S(\gamma_{[a, t]})=S(\gamma_{[a,
t_0]})+l_{\gamma}(t).$ Крім того, нижче ми використовуємо позначення
$l_{\gamma}(t)$ замість $s_{\gamma}(t),$ якщо непорозуміння
неможливе, $l_{\gamma}(t)=S\left(\gamma, [a,t]\right),$ де
$S(\gamma, [a,t])$ -- довжина кривої $\gamma|_{[a, t]}.$
\end{remark}

\medskip
Нехай $\alpha:[a,b]\rightarrow {\Bbb R}^n$ -- спрямлювана замкнена
крива в ${\Bbb R}^n,$ $n\geqslant 2,$ $l(\alpha)$ -- її довжина.
{\it Нормальним представленням кривої $\alpha$ } називається крива
$\alpha^0:[0, l(\alpha)]\rightarrow {\Bbb R}^n,$ така що
$\alpha(t)=\alpha^0\left(S\left(\alpha, [a, t]\right)\right)=
\alpha^0\circ l_{\alpha}(t).$ Відзначимо, що така крива
$\alpha^{\,0}$ існує та єдина, і що в цьому випадку
$S\left(\alpha^0, [0, t]\right)=t$ при $t\in [0, l(\alpha)],$ див.
теорему~2.4 \cite{Va}.

Наступне означення може бути знайденим у \cite[2.5 п. 2 розд.
I]{Va}. Нехай $\alpha:[a,b]\rightarrow {\Bbb R}^n$ -- спрямлювана
замкнена крива в ${\Bbb R}^n,$ $n\geqslant 2.$ Відображення
$f:|\alpha|\rightarrow {\Bbb R}^n$ називається {\it абсолютно
неперервним на $\alpha,$} якщо композиція $f\circ\alpha^0$ є
абсолютно неперервною на інтервалі $[0, l(\alpha)],$ де $l(\alpha)$
позначає довжину кривої $\alpha,$ а $\alpha^0$ -- її нормальне
представлення.

\medskip
\begin{remark}\label{rem5.1A}
Зауважимо, що абсолютна неперервність відображення $f$ на локально
спрямлюваній кривій $\gamma$ тягне локальну спрямлюваність кривої
$\gamma^{\,\prime}=f\circ\gamma,$ див. \cite[теореми 2.6 і 5.3]{Va}.
\end{remark}

Нехай $\alpha$ і $\beta$ -- криві в ${\Bbb R}^n,$ тоді запис
$\alpha\subset\beta$ позначає, що $\alpha$ є підкривою кривої
$\beta.$ Надалі $I$ позначає відкритий, замкнений, або
напіввідкритий інтервал числової осі. Наступне означення див.,
напр., у \cite[п. 5 розд. II]{Ri}. Нехай $f:D\rightarrow {\Bbb R}^n$
слабко нульвимірне відображення, $\beta:I_0\rightarrow {\Bbb R}^n$
замкнена спрямлювана крива і $\alpha:I\rightarrow D$ крива така, що
$f\circ \alpha\subset \beta.$ Якщо функція довжини
$l_{\beta}:I_0\rightarrow [0, l(\beta)]$ є сталою на певному
інтервалі $J\subset I,$ то $\beta$ є сталою на $J$ і, в силу слабкої
нульвимірності відображення $f,$ крива $\alpha$ також є сталою на
$J.$ Отже, існує єдина функція $\alpha^{\,*}:l_\beta(I)\rightarrow
D$ така, що $\alpha=\alpha^{\,*}\circ (l_\beta|_I).$ Будемо
говорити, що $\alpha^{\,*}$ є {\it $f$-пред\-став\-лен\-ням кривої
$\alpha$ відносно $\beta.$ }

Нехай $X$ та $Y$ -- два простори з мірами $\mu$ і $\mu^{\,\prime},$
відповідно. Будемо говорити, що відображення $f:X\rightarrow Y$ має
{\it $N$-властивість Лузіна,} якщо з умови $\mu(E)=0$ випливає, що
$\mu^{\,\prime}(f(E))=0.$ Аналогічно, будемо говорити, що
відображення $f:X\rightarrow Y$ має {\it $N^{-1}$-властивість
Лузіна,} якщо умова $\mu^{\,\prime}(E)=0$ тягне, що
$\mu(f^{\,-1}(E))=0.$ Покладемо в точках диференційовності $x\in D$
відображення $f$
$$l(f^{\,\prime}(x))\,=\,\min\limits_{h\in {\Bbb R}^n
\backslash \{0\}} \frac {|f^{\,\prime}(x)h|}{|h|}\,,$$
\begin{equation}\label{eq5_a}
\Vert f^{\,\prime}(x)\Vert\,=\,\max\limits_{h\in {\Bbb R}^n
\backslash \{0\}} \frac {|f^{\,\prime}(x)h|}{|h|}\,,
\end{equation}
$$J(x,f)=\det
f^{\,\prime}(x)\,.$$

\medskip
Надалі ми говоримо, що деяка властивість $P$ виконується для {\it
$p$-майже всіх кривих в області $D$}, якщо ця властивість може
порушуватись лише для деякої сім'ї кривих $\Gamma_0$ у $D$ такої, що
$M_p(\Gamma_0)=0,$ де а $M_p(\Gamma_0)$ позначає $p$-модуль сім'ї
кривих $\Gamma_0$ (див.~\cite[розд.~6]{Va}).  Будемо говорити, що
ві\-доб\-ра\-жен\-ня $f:D\rightarrow {\Bbb R}^n$ має {\it
$ACP$--влас\-ти\-вість відносно $p$-модуля,} пишемо $f\in ACP_p,$
якщо функція довжини $L_{\gamma, f}$ є абсолютно неперервною на всіх
замкнених інтервалах $\Delta_{\gamma}$ для $p$-майже всіх кривих
$\gamma$ у $D.$ Іншими словами, $f\in ACP_p\Leftrightarrow$ звуження
$f|_{\gamma}$ є локально абсолютно неперервною функцією для майже
всіх кривих $\gamma.$

\medskip
Нехай $y_0\in {\Bbb R}^n,$ $0<r_1<r_2<\infty$ і
\begin{equation}\label{eq1**}
A=A(y_0, r_1,r_2)=\left\{ y\,\in\,{\Bbb R}^n:
r_1<|y-y_0|<r_2\right\}\,.\end{equation}
Для заданих множин $E,$ $F\subset\overline{{\Bbb R}^n}$ і області
$D\subset {\Bbb R}^n$ позначимо через $\Gamma(E,F,D)$ сім'ю всіх
кривих $\gamma:[a,b]\rightarrow \overline{{\Bbb R}^n}$ таких, що
$\gamma(a)\in E,\gamma(b)\in\,F$ і $\gamma(t)\in D$ при $t \in [a,
b].$ Якщо $f:D\rightarrow {\Bbb R}^n$ -- задане відображення,
$y_0\in \overline{f(D)}\setminus\{\infty\},$ і
$0<r_1<r_2<r_0=\sup\limits_{y\in f(D)}|y-y_0|,$ то через
$\Gamma_f(y_0, r_1, r_2)$ позначимо сім'ю всіх кривих $\gamma$ в
області $D$ таких, що $f(\gamma)\in \Gamma(S(y_0, r_1), S(y_0, r_2),
A(y_0,r_1,r_2)).$ Нехай $Q_*:{\Bbb R}^n\rightarrow [0, \infty]$ --
вимірна за Лебегом функція. Будемо говорити, що {\it $f$ задовольняє
обернену нерівність Полецького в точці $y_0\in
\overline{f(D)}\setminus\{\infty\}$ відносно $p$-модуля,} якщо
співвідношення
\begin{equation}\label{eq2*A}
M_p(\Gamma_f(y_0, r_1, r_2))\leqslant
\int\limits_{A(y_0,r_1,r_2)\cap f(D)} Q_*(y)\cdot \eta^{p}
(|y-y_0|)\, dm(y)
\end{equation}
виконується для довільної вимірної за Лебегом функції $\eta:
(r_1,r_2)\rightarrow [0,\infty ]$ такій, що
\begin{equation}\label{eqA2}
\int\limits_{r_1}^{r_2}\eta(r)\, dr\geqslant 1\,.
\end{equation}
Справедлива наступна

\begin{theorem}\label{th1} {\sl Нехай $p>1,$ $f:D\rightarrow {\Bbb R}^n$ --
диференційовне майже скрізь відображення, яке має $N$- та
$N^{\,-1}$-властивості Лузіна відносно лебегової міри в ${\Bbb
R}^n,$ причому $f\in ACP_p(D).$ Нехай $y_0\in
\overline{f(D)}\setminus \{\infty\}.$ Покладемо
\begin{equation}\label{eq1}
K_{CT, p, y_0}(y, f)=\sum\limits_{x\in
f^{\,-1}(y)}\frac{\left(\sup\limits_{|h|=1}\left|\left(f^{\,\prime}(x)h,
\frac{f(x)-y_0}{|f(x)-y_0|}\right)\right|\right)^p}{|J(x, f)|}\,.
\end{equation}
Тоді відображення $f$ задовольняє обернену нерівність
Полецького~(\ref{eq2*A}) в точці $y_0$ при $Q_*(y):=K_{CT, p,
y_0}(y, f).$
 }
\end{theorem}
Тут і надалі $(x, y)$ позначає скалярний добуток векторів $x, y\in
{\Bbb R}^n.$ Величину $K_{CT, p, y_0}(y, f)$ в~(\ref{eq1}) назвемо
{\it котангенсальною дилатацію відображення $f$ порядку $p$ в точці
$y_0$}. Докладніше про походження цього терміну буде згадано нижче
за текстом.

\section{Доведення теореми~\ref{th1}}

У значній мірі доведення теореми~\ref{th1} спирається на підхід,
використаний при доведенні~\cite[теорема~2.1]{SalSev$_1$}, див.
також \cite{MRSY}, або \cite[теорема~6.1]{MRSY$_1$}. Перед тим, як
безпосередньо переходити до нього, наведемо декілька корисних
тверджень.

\medskip
Нехай $E$ -- множина в ${\Bbb R}^n$ і $\gamma :\Delta\rightarrow
{\Bbb R}^n$ деяка крива. Позначимо через $\gamma\cap
E\,=\,\gamma(\Delta)\cap E.$ Нехай крива $\gamma$ є локально
спрямлюваною і функцію довжини $l_{\gamma}(t)$ визначено вище.
Покладемо
$$
l(\gamma\cap E):= {\rm mes}_1\,(E_ {\gamma}), \quad E_ {\gamma} =
l_{\gamma}(\gamma ^{\,-1}(E))\,.
$$
Тут, як і всюди вище, ${\rm mes}_1\,(A)$ позначає довжину (лінійну
міру Лебега) множини $A\subset {\Bbb R}.$ Зауважимо, що
$$E_ {\gamma} = \gamma_0^{\,-1}\left(E\right)\,,$$
де $\gamma _0 :\Delta _{\gamma}\rightarrow {\Bbb R}^n$ -- натуральна
параметризація кривої $\gamma,$  і що
$$l\left(\gamma\cap E\right) = \int\limits_{\gamma} \chi_E(x)\,|dx|
= \int\limits_{\Delta _{\gamma}} \chi _{E_\gamma }(s)\,dm_1(s)\,.$$

\medskip
Наступне твердження може бути знайдено в~\cite[теорема~9.1,
$k=~1$]{MRSY}.

\medskip
\begin{proposition}\label{pr2}{\sl\,
Нехай $E$ -- підмножина області $D\subset{\Bbb R}^n,$ $n\geqslant
2,$ $p\geqslant 1.$ Тоді множина $E$ є вимірною за Лебегом тоді і
тільки тоді, коли множина $\gamma\cap E$ є вимірною для $p$-майже
всіх кривих $\gamma$ в $D.$ Більше того, $m(E)=0$ тоді і тільки
тоді, коли
$$l(\gamma\cap E)=0$$
для $p$-майже всіх кривих $\gamma$ в $D.$}
\end{proposition}

\medskip
Наступне твердження дає характеристику відображень зі скінченим
спотворенням довжини на мові абсолютної неперервності кривих і
повністю доведено в~\cite[Пропозиція~2.1]{SalSev$_1$}.

\medskip
\begin{proposition}\label{pr11}
{\sl\, Відображення $f:D\rightarrow {\Bbb R}^n$ має
$ACP_p$-влас\-ти\-вість тоді і тільки тоді, коли $f$ є абсолютно
неперервним на майже всіх замкнених кривих (тобто $f\circ\gamma^0$ є
спрямлюваною і абсолютно неперервною для майже всіх замкнених кривих
$\gamma$).}
\end{proposition}

\medskip
Нагадаємо, що ві\-доб\-ра\-жен\-ня $\varphi:X\rightarrow Y$ між
метричними просторами $X$ і $Y$ називається {\it ліпшицевим,} якщо
${\rm dist}\,\left(\varphi(x_1),\,\varphi(x_2)\right)\leqslant
M\cdot {\rm dist}\,\left(x_1,\,x_2\right)$
для певної сталої $M<\infty$ і всіх $x_1, x_2 \in X.$ Говорять, що
ві\-доб\-ра\-жен\-ня $\varphi:X\rightarrow Y$ {\it є біліпшицевим,}
якщо: по-перше, воно є ліпшицевим, по-друге,
$M^{\,*}\cdot {\rm dist}\,\left(x_1,\,x_2\right)\leqslant{\rm
dist}\,\left(\varphi(x_1),\,\varphi(x_2)\right)$
для певної сталої $M^{\,*}>0$ і всіх $x_1, x_2 \in\,X.$ Надалі
$X:=D$ -- область в ${\Bbb R}^n,$ $Y:={\Bbb R}^n$ і ${\rm
dist\,}(x_1, x_2):=|x_1-x_2|.$

\medskip
Наступний результат, що буде використовуватися у подальшому, було
отримано в \cite[лема~3.20, наслідок~3.14]{MRSY$_1$}, див. також
\cite[леми 8.2 і 8.3, наслідок~8.1]{MRSY}.

\medskip
\begin{proposition}\label{lem2.4}
{\sl\, Нехай $f:D\rightarrow {\Bbb R}^n$ -- ві\-доб\-ра\-жен\-ня,
яке є диференційовним майже скрізь та має $N$- та
$N^{\,-1}$-властивості. Тоді існує зчисленна послідовність
компактних множин $C_k^{\,*}\subset D,$ така що $m(B)=0,$ де
$B=D\setminus \bigcup\limits_{k=1}^{\infty} C_k^{\,*}$ і
$f|_{C_k^{\,*}}$ є взаємно однозначним та біліпшицевим для кожного
$k=1,2,\ldots .$ Більше того, $f$ диференційовне для всіх $x\in
C_k^{\,*},$ причому $J(x,f)\ne 0$ на $C_k^{\,*}.$}
\end{proposition}

\medskip
{\it Доведення теореми~\ref{th1}}. Припустимо, що $B_0$ і $C_k^*,$
$k=1,2,\ldots ,$ є такими, як у твердженні~\ref{lem2.4}. Покладаючи
по індукції $B_1=C_1^*,$ $B_2=C_2^*\setminus B_1,\ldots ,$ і
\begin{equation} \label{eq7.3.7y} B_k=C_k^*\setminus
\bigcup\limits_{l=1}\limits^{k-1}B_l\,,
\end{equation}
ми отримаємо зчисленне покриття області $D,$ яке складається з
попарно непересічних борелевих множин $B_k, k=0,1,2,\ldots ,$ таких
що $m(B_0)=0,$ $B_0=D\setminus \bigcup\limits_{k=1}^{\infty} B_k.$
Зауважимо, що $\gamma^{\,0}(s)\not\in B_0$ для майже всіх $s$ і
$p$-майже всіх замкнених кривих $\gamma\in \Gamma,$ див.
твердження~\ref{pr2}; тут, як звично, $\gamma^{0}(s)$ позначає
нормальне зображення кривої $\gamma.$ За твердженням~\ref{pr11}
крива $f\circ\gamma^{0}$ є спрямлюваною і абсолютно неперервною для
$p$-майже всіх кривих $\gamma\in\Gamma.$

\medskip
Зафіксуємо тепер $0<r_1<r_2<\infty$ і вимірну за Лебегом функцію
$\eta:(r_1, r_2)\rightarrow [0, \infty]$ таку, що
$\int\limits_{r_1}^{r_2}\eta(t)\,dt\geqslant 1.$ Ми можемо вважати
функцію $\rho$ борелевою, бо за теоремою Лузіна існує борелева
функція $\eta_1,$ яка дорівнює функції $\eta$ майже скрізь (див.,
напр., \cite[розд.~2.3.6]{Fe}). Крім того, згідно зроблених вище
зауважень, ми можемо вважати, що кожна крива $f\circ \gamma,$
$\gamma\in \Gamma_f(y_0, r_1, r_2),$ є спрямлюваною і абсолютно
неперервною. Ми також можемо вважати спрямлюваними всі криві
$\gamma$ самої сім'ї $\Gamma_f(y_0, r_1, r_2).$

\medskip
Покладемо тепер $\rho(x)=0,$ якщо $|f(x)-y_0|<r_1,$ або
$|f(x)-y_0|>r_2;$ крім того, при $r_1<|f(x)-y_0|<r_2$ покладемо
\begin{equation}\label{eq2}
\rho (x)= \left \{\begin{array}{rr}
\eta(|f(x)-y_0|)\sup\limits_{|h|=1}\left|\left(f^{\,\prime}(x)h,
\frac{f(x)-y_0}{|f(x)-y_0|}\right)\right|, &  x\in D\setminus B_0, \\
0, & x\not\in D\setminus B_0.\end{array} \right.
\end{equation}
Можна показати, що функція $\rho$ -- борелева. Тоді в силу зроблених
вище зауважень будемо мати, що
$$\int\limits_{\gamma} \rho(x) |dx|=$$
\begin{equation}\label{eq3}
=\int\limits_{0}^{l(\gamma)}
\eta(|f(\gamma^0(s))-y_0|)\sup\limits_{|h|=1}\left|\left(f^{\,\prime}(\gamma^0(s))h,
\frac{f(\gamma^0(s))-y_0}{|f(\gamma^0(s))-y_0|}\right)\right|\,ds\,,
\end{equation}
де $\gamma^0,$ як звично, позначає натуральне зображення кривої
$\gamma.$ Нехай $s_1$ -- натуральний параметр на кривій $f\circ
\gamma.$ Тоді $s=s(s_1),$ причому за означенням класу $ACP_p$
функція $s_1=s_1(s)$ є абсолютно неперервною для $p$-майже всіх
кривих сім'ї $\Gamma_f(y_0, r_1, r_2).$ В такому випадку, за
теоремою Пономарьова обернена функція $s=s(s_1)$ має нерівну нулю
похідну майже скрізь, див. \cite[теорема~1]{Pon}. Зауважимо також,
що $\frac{ds_1}{ds}=|\left(f(\gamma^0(s))\right)^{\prime}_s|$ майже
скрізь, див.~\cite[теорема~1.3]{Va}. В свою чергу, за теоремою про
похідну складеної функції
\begin{equation}\label{eq5}
\frac{ds_1}{ds}=|\left(f(\gamma^0(s))\right)^{\prime}_s|=
|f({\gamma^0}^{\,\prime}(s))\cdot {\gamma^0}^{\,\prime}(s)|\,,
\end{equation}
причому $|{\gamma^0}^{\,\prime}(s)|=1$ (бо за~\cite[теорема~1.3]{Va}
$|{\gamma^0}^{\,\prime}(s)|=\frac{ds}{ds}=1$ майже скрізь). З
рівності~(\ref{eq5}) випливає, що
\begin{equation}\label{eq6}
\frac{ds_1}{ds}\geqslant l(f^{\,\prime}(\gamma^0(s)))>0
\end{equation}
при майже всіх $s,$ де $(f^{\,\prime}(\gamma^0(s)))$ визначено
першою рівністю у~(\ref{eq5_a}), а $l(f^{\,\prime}(\gamma^0(s)))>0$
з огляду на те, що $\gamma^{\,0}(s)\not\in B_0$ для майже всіх $s.$
Враховуючи~\cite[теорема~3.2.6]{Fe} та рівності у~(\ref{eq5}),
перетворимо вираз у~(\ref{eq3}) наступним чином:
$$\int\limits_{0}^{l(\gamma)}
\eta(|f(\gamma^0(s))-y_0|)\sup\limits_{|h|=1}\left|\left(f^{\,\prime}(\gamma^0(s))h,
\frac{f(\gamma^0(s))-y_0}{|f(\gamma^0(s))-y_0|}\right)\right|\,ds=$$$$=
\int\limits_{0}^{l(\gamma)}
\eta(|f(\gamma^0(s))-y_0|)\sup\limits_{|h|=1}\left|\left(f^{\,\prime}(\gamma^0(s))h,
\frac{f(\gamma^0(s))-y_0}{|f(\gamma^0(s))-y_0|}\right)\right|\cdot
\frac{\frac{ds_1}{ds}}{\frac{ds_1}{ds}}\,ds=$$
\begin{equation}\label{eq7}=\int\limits_{0}^{l(f(\gamma))}
\eta(r)\frac{\sup\limits_{|h|=1}\left|\left(f^{\,\prime}(\gamma^0(s(s_1)))h,
\frac{f(\gamma^0(s(s_1)))-y_0}{|f(\gamma^0(s(s_1)))-y_0|}\right)\right|}
{\frac{ds_1}{ds}}\,ds_1\geqslant
\end{equation}
$$\geqslant\int\limits_{0}^{l(f(\gamma))}
\eta(r)\frac{\left|\left(f^{\,\prime}(\gamma^0(s(s_1))){\gamma^0}^{\,\prime}(s(s_1)),
\frac{f(\gamma^0(s(s_1)))-y_0}{|f(\gamma^0(s(s_1)))-y_0|}\right)\right|}
{\frac{ds_1}{ds}}\,ds_1\,,$$
де $r=r(s_1)=|f(\gamma^0(s(s_1)))-y_0|.$ Зробимо заміну змінних
$r=r(s_1).$ Тоді
$$\frac{dr}{ds_1}=\frac{dr}{ds}\cdot \frac{ds}{ds_1}=$$
\begin{equation}\label{eq8}
= \left(f^{\,\prime}(\gamma^0(s(s_1))){\gamma^0}^{\,\prime}(s(s_1)),
\frac{f(\gamma^0(s(s_1)))-y_0}{|f(\gamma^0(s(s_1)))-y_0|}\right)\cdot
\frac{ds}{ds_1}\,.
\end{equation}
Оскільки $\frac{ds}{ds_1}\ne 0$ майже скрізь, можна застосувати
теоремою про похідну оберненої функції. Згідно цієї теореми
$\frac{ds}{ds_1}=(\frac{ds_1}{ds})^{\,-1}.$ Тоді з~(\ref{eq8})
випливає, що
\begin{equation}\label{eq9}
\frac{dr}{ds_1}=
\left(f^{\,\prime}(\gamma^0(s(s_1))){\gamma^0}^{\,\prime}(s(s_1)),
\frac{f(\gamma^0(s(s_1)))-y_0}{|f(\gamma^0(s(s_1)))-y_0|}\right)
\cdot\frac{1}{\frac{ds_1}{ds}}\,.
\end{equation}
Поєднуючи тепер~(\ref{eq7}), (\ref{eq8}) та~(\ref{eq9}), будемо
мати:
$$\int\limits_{\gamma} \rho(x) |dx|\geqslant$$
$$\geqslant\int\limits_{0}^{l(\gamma)}
\eta(|f(\gamma^0(s))-y_0|)\sup\limits_{|h|=1}\left|\left(f^{\,\prime}(\gamma^0(s))h,
\frac{f(\gamma^0(s))-y_0}{|f(\gamma^0(s))-y_0|}\right)\right|\,ds\,\geqslant$$
\begin{equation}\label{eq10}
\geqslant\int\limits_{0}^{l(f(\gamma))} \eta(r(s_1))\,dr(s_1)\,.
\end{equation}
Зауважимо, що функція $r=r(s_1)$ є абсолютно неперервною по змінній
$s_1.$ Дійсно,
$r=r(s_1)=|f(\gamma^0(s(s_1)))-y_0|=|(f(\gamma))^0(s_1)-y_0|,$
причому крива $(f(\gamma))^0(s_1)-y_0$ абсолютно неперервна по $s_1$
як крива, що є абсолютно неперервною відносно свого натурального
параметру, а функція $\varphi(x)=|x|$ є ліпшицевою. Тоді
за~\cite[теорема~3.2.6]{Fe} будемо мати:
\begin{equation}\label{eq11}
\int\limits_{0}^{l(f(\gamma))}
\eta(r(s_1))\,dr(s_1)=\int\limits_{r_1}^{r_2}\eta(t)\,dt\geqslant
1\,.
\end{equation}
Поєднуючи~(\ref{eq10}) і~(\ref{eq11}), ми отримаємо, що
\begin{equation}\label{eq12}
\int\limits_{\gamma} \rho(x) |dx|\geqslant 1
\end{equation}
для $p$-майже всіх $\gamma\in \Gamma_f(y_0, r_1, r_2).$ Отже,
$\rho\in {\rm adm\,}\Gamma_f(y_0, r_1, r_2)\setminus \Gamma_0,$ де
$M_p(\Gamma_0)=0.$ (Тут і надалі $\rho\in {\rm adm\,}\Gamma,$ якщо
нерівність~(\ref{eq12}) виконується для довільної локально
спрямлюваної кривої $\gamma\in\Gamma$). Отже,
\begin{equation} \label{eq7.3.7z} M_p(\Gamma_f(y_0, r_1, r_2))\leqslant\int\limits_{D}
\rho^p(x) dm(x). \end{equation}
Зауважимо, що $\rho =\sum\limits_{k=1}\limits^{\infty}\rho_k$, де
функції $\rho_k = \rho\cdot\chi_{B_k}$ мають попарно непересічні
носії. Позначимо
$$K^*_{T}(x, y_0):=\frac{\left(\sup\limits_{|h|=1}\left|\left(f^{\,\prime}(x)h,
\frac{f(x)-y_0}{|f(x)-y_0|}\right)\right|\right)^p}{|J(x, f)|}\,.$$
Тоді за~\cite[теорема~3.2.5 при $m=n$]{Fe}
$$\int\limits_{f(B_k)\cap A(y_0, r_1, r_2)} K^*_{T}(f_k^{-1}(y), y_0)
\cdot\eta^p(|y-y_0|)\,dm(y)= $$
$$=\int\limits_{B_k} K^*_{T}(x, y_0)
\cdot\eta^p(|f(x)-y_0|)|J(x, f)|dm(x)=$$
\begin{equation}\label{eq7.3.7x}
=\int\limits_{B_k}
\left(\sup\limits_{|h|=1}\left|\left(f^{\,\prime}(x)h,
\frac{f(x)-y_0}{|f(x)-y_0|}\right)\right|\right)^p
\eta^p(|f(x)-y_0|) \,dm(x)=
\end{equation}
$$=\int\limits_{D}\rho_k^p(x)\,dm(x)\,,$$
де кожне з відображень $f_k=f|_{B_k},$ $k=1,2,\ldots $ є ін'єктивним
за побудовою. Остаточно, оскільки
$$K_{CT, p, y_0}(y, f)=\sum\limits_{x\in
f^{\,-1}(y)}\frac{\left(\sup\limits_{|h|=1}\left|\left(f^{\,\prime}(x)h,
\frac{f(x)-y_0}{|f(x)-y_0|}\right)\right|\right)^p}{|J(x,
f)|}=\sum\limits_{x\in f^{\,-1}(y)}K^*_{T}(x, y_0)\,,$$
то сумуючи по $k=1,2,\ldots $ у (\ref{eq7.3.7x}) і застосовуючи
теорему Лебега про збіжність позитивних рядів,
див.~\cite[теорема~I.12.3]{Sa}, ми отримаємо, що
$$\int\limits_{A(y_0, r_1, r_2)\cap f(D)} K_{CT, p, y_0}(y, f)
\cdot\eta^p(|y-y_0|)\,dm(y)=$$$$=
\sum\limits_{k=1}\limits^{\infty}\int\limits_{D}\rho_k^p(x)dm(x)
\geqslant M_p(\Gamma_f(y_0, r_1, r_2))\,.\,\Box$$

\medskip
\begin{remark}\label{rem1}
З огляду на нерівність Коші-Буняковського,
\begin{equation}\label{eq13}
K_{CT, p, y_0}(y, f)=\sum\limits_{x\in
f^{\,-1}(y)}\frac{\left(\sup\limits_{|h|=1}\left|\left(f^{\,\prime}(x)h,
\frac{f(x)-y_0}{|f(x)-y_0|}\right)\right|\right)^p}{|J(x,
f)|}\leqslant$$$$\leqslant \sum\limits_{x\in f^{\,-1}(y)}\frac{\Vert
f^{\,\prime}(x)\Vert^p}{|J(x, f)|}:=K_{I,p}(y, f^{\,-1})\,.
\end{equation}
Тому нерівність
$$M_p(\Gamma_f(y_0, r_1, r_2))\leqslant$$
\begin{equation}\label{eq14}
\leqslant \int\limits_{A(y_0,r_1,r_2)\cap f(D)} K_{CT, p, y_0}(y,
f)\cdot \eta^{p} (|y-y_0|)\, dm(y)\,,
\end{equation}
яка отримана в теоремі~\ref{th1}, є більш сильною у порівнянні з
нерівністю
$$M_p(\Gamma_f(y_0, r_1, r_2))\leqslant$$
\begin{equation}\label{eq15} \leqslant
\int\limits_{A(y_0,r_1,r_2)\cap f(D)} K_{I,p}(y, f^{\,-1})\cdot
\eta^{p} (|y-y_0|)\, dm(y)\,,\end{equation}
яку було отримано раніше в деяких роботах (див., напр.,
\cite[теорема~8.5]{MRSY}, \cite[теорема~2.1]{SalSev$_1$}). Слід
розуміти, що у буквальному сенсі нерівність~(\ref{eq14}) не є більш
загальною у порівнянні з результатами згаданих робіт, бо в останніх
нерівність~(\ref{eq15}) сформульована в більш абстрактному вигляді:
\begin{equation}\label{eq16} M_p(\Gamma)\leqslant
\int\limits_{f(E)} K_{I,p}(y, f^{\,-1}, E)\cdot \rho_*^{p}(y)\,
dm(y)\,,\end{equation}
де сім'я $\Gamma$ є довільною і належить до довільної вимірної
множини $E\subset D,$ а $\rho_*\in {\rm adm}f(\Gamma).$ Отже, сім'я
кривих $\Gamma$ в~(\ref{eq16}) -- довільна, а в~(\ref{eq14}) --
<<спеціального вигляду>>, яка чітко прив'язана до точки $y_0.$ В той
самий час, як було сказано, для цієї <<спеціальної>> сім'ї кривих
нерівність~(\ref{eq14}) сильніша за~(\ref{eq15}), або, що те саме,
за~(\ref{eq16}), де функції $\rho_*(y)$ також мають окремий вигляд
$\eta(|y-y_0|).$
\end{remark}

\medskip
\begin{remark}\label{rem2}
В роботі~\cite{RSY}, див. також \cite{GG} і \cite{SalSev},
впроваджено означення {\it тангенсальної (дотичної) дилатації}. Так
називається наступна величина:
\begin{equation}\label{equa11}
D_f(x, x_0)=\frac{|J(x ,f)|}{l_f^n(x, x_0)}\,,
\end{equation}
де $$l_f(x, x_0)=\min\limits_{|h|=1}\frac{|\partial_h
f(x)|}{\left|\left(h, \frac{x-x_0}{|x-x_0|}\right)\right|}\,,$$
$\partial_h f(x)=\lim\limits_{t\rightarrow
+0}\frac{f(x+th)-f(x)}{t}.$ Покажемо, що якщо відображення $f$ є
гомеоморфізмом, то в точках $x\in D,$ в яких $f$ є невироджено
диференційовним
\begin{equation}\label{eq17}
K_{CT, n, x_0}(f(x), f^{\,-1})=D_f(x, x_0)\,,\quad y_0=f(x_0)\,.
\end{equation}
Дійсно,
$$K_{CT, n, x_0}(f(x), f^{\,-1})=
\frac{\left(\sup\limits_{|h|=1}\left|\left({f^{\,-1}}^{\,\prime}(f(x))h,
\frac{f^{\,-1}(f(x))-x_0}{|f^{\,-1}(f(x))-x_0|}\right)\right|\right)^n}{|J(f(x),
f^{\,-1})|}=$$
$$=\left(\sup\limits_{|h|=1}\left|\left((f^{\,\prime}(x))^{\,-1}h,
\frac{x-x_0}{|x-x_0|}\right)\right|\right)^n\cdot |J(x, f)|=$$
$$=\left(\sup\limits_{h\in {\Bbb R}^n\setminus\{0\}}\left|
\left((f^{\,\prime}(x))^{\,-1}\frac{h}{|h|},
\frac{x-x_0}{|x-x_0|}\right)\right|\right)^n\cdot |J(x, f)|=$$
$$=\left(\sup\limits_{h\in {\Bbb R}^n\setminus\{0\}}\frac{\left|
\left((f^{\,\prime}(x))^{\,-1}h,
\frac{x-x_0}{|x-x_0|}\right)\right|}{|h|}\right)^n\cdot |J(x, f)|=$$
$$=\left(\sup\limits_{H\in {\Bbb R}^n\setminus\{0\}}\left|\frac{
\left(H,
\frac{x-x_0}{|x-x_0|}\right)}{|f^{\,\prime}(x)H|}\right|\right)^n\cdot
|J(x, f)|=\left(\sup\limits_{H\in {\Bbb
R}^n\setminus\{0\}}\frac{\left|\left(\frac{H}{|H|},
\frac{x-x_0}{|x-x_0|}\right)\right|}{|f^{\,\prime}(x)\frac{H}{|H|}|}\right)^n\cdot
|J(x, f)|=$$
$$=\left(\sup\limits_{|H|=1}\frac{\left|\left(H,
\frac{x-x_0}{|x-x_0|}\right)\right|}{|f^{\,\prime}(x)H|}\right)^n\cdot
|J(x, f)|=
\frac{1}{\left(\inf\limits_{|H|=1}\frac{1}{\left(\frac{\left|\left(H,
\frac{x-x_0}{|x-x_0|}\right)\right|}{|f^{\,\prime}(x)H|}\right)}\right)^n}\cdot
|J(x, f)|=$$$$=\frac{|J(x,
f)|}{\left(\inf\limits_{|H|=1}\frac{|f^{\,\prime}(x)H|}{\left|\left(H,
\frac{x-x_0}{|x-x_0|}\right)\right|}\right)^n}=D_f(x, x_0)\,.$$
Отже, співвідношення~(\ref{eq17}) встановлено.
\end{remark}

\medskip
\begin{remark}\label{rem3}
Нехай відображення $f$ -- гомеоморфізм. Як вже було зазначено в
зауваженні~\ref{rem1},
\begin{equation}\label{eq19}
K_{CT, p, y_0}(y, f)\leqslant K_{I,p}(y, f^{\,-1})\,.
\end{equation}
Зауважимо, що остання нерівність, взагалі кажучи, є строгою на
деякій множині додатної міри. Дійсно, нехай, наприклад, $p=n=2.$
Тоді згідно зауваження~\ref{rem2}
\begin{equation}\label{eq18}
K_{CT, 2, x_0}(f(x), f^{\,-1})=D_f(x, x_0)\,,\quad y_0=f(x_0)\,,
x_0\in D\subset{\Bbb C}\,.
\end{equation}
Зауважимо, що
\begin{equation}\label{eq20}
D_f(x,
x_0)=\frac{\left|1-\frac{\overline{x-x_0}}{x-x_0}\mu(x)\right|^2}{1-|\mu(x)|^2}\,,
\end{equation}
де $\mu(x)=\frac{f_{\overline x}}{f_x},$ $f_x\ne 0,$ $\mu(x)=0$ при
$f_x=0,$ $f_x=(f_{x_1}-if_{x_2})/2,$
$f_{\overline{x}}=(f_{x_1}+if_{x_2})/2,$ $x=x_1=ix_2,$ $i^2=-1$
(див. \cite[лема~11.2]{MRSY}). Нехай, наприклад, $x_0=0,$ тоді
\begin{equation}\label{eq21}
D_f(x,
0)=\frac{\left|1-\frac{\overline{x}}{x}\mu(x)\right|^2}{1-|\mu(x)|^2}\,.
\end{equation}
Якщо, наприклад, $\mu(x)=\frac{1}{2}\cdot \frac{x}{\overline{x}},$
то з~(\ref{eq21}) отримаємо, що
\begin{equation}\label{eq22}
D_f(x,
0)=\frac{\left|1-\frac{\overline{x}}{x}\mu(x)\right|^2}{1-|\mu(x)|^2}=\frac{\frac
14}{\frac 34}=\frac 13<1\,.
\end{equation}
Тоді з огляду на~(\ref{eq18}) $K_{CT, 2, x_0}(f(x),
f^{\,-1})=\frac{1}{3},$ проте це робить нерівність~(\ref{eq19})
строгою, бо права частина~(\ref{eq19}) при $p=n$ завжди не менша
одиниці (див., напр.,\cite[співвідношення (4.8)--(4.10),
гл.~I]{Re}).
\end{remark}

КОНТАКТНА ІНФОРМАЦІЯ

\medskip
\noindent{{\bf Євген Олександрович Севостьянов} \\
{\bf 1.} Житомирський державний університет ім.\ І.~Франко\\
кафедра математичного аналізу, вул. Велика Бердичівська, 40 \\
м.~Житомир, Україна, 10 008 \\
{\bf 2.} Інститут прикладної математики і механіки
НАН України, \\
вул.~Добровольського, 1 \\
м.~Слов'янськ, Україна, 84 100\\
e-mail: esevostyanov2009@gmail.com}

\medskip
\noindent{{\bf Валерій Андрійович Таргонський} \\
Житомирський державний університет ім.\ І.~Франко\\
кафедра математичного аналізу, вул. Велика Бердичівська, 40 \\
м.~Житомир, Україна, 10 008 \\
e-mail: w.targonsk@gmail.com }

\end{document}